\documentstyle{article}
\begin{document}

\centerline {\bf IDENTICAL AND NONIDENTICAL RELATIONS.}
\centerline {\bf NONDEGENERATE AND DEGENERATE}
\centerline {\bf TRANSFORMATIONS}
\bigskip
\centerline {{\bf (Properties of skew-symmetric differential forms)}}
\bigskip
\centerline {L.I. Petrova}
\bigskip
\centerline{{\it Moscow State University, Russia, e-mail: ptr@cs.msu.su}}
\bigskip

Identical relations occur in various branches of mathematics and 
mathematical physics. The Cauchy-Riemann relations, 
characteristical and canonical relations, the Bianchi identities 
and others are examples of identical relations. It can be shown 
that all these relations express either the conditions of closure 
of exterior (skew-symmetric) differential forms and corresponding 
dual forms or the properties of closed exterior forms. Since the 
closed differential forms are invariant under all transformations, 
which conserve the differential (these are gauge transformations: 
unitary, canonical, gradient and others), from this it follows 
that identical relations are a mathematical representation of 
relevant invariant and covariant objects, which are of great 
functional and utilitarian importance. 

The theory of exterior differential forms, which describes invariant 
objects using identical relations, cannot answer the question of how 
do invariant objects appear and what does these objects generate? 

The answer to this question can be obtained using the skew-symmetric 
differential forms, which possess the evolutionary properties. The 
mathematical apparatus of such evolutionary forms contains nonidentical 
relations, from which the identical relations corresponding to invariant 
objects are obtained with the help of degenerate transformations. 

Due to such potentialities, the mathematical apparatus of 
skew-symmetrical differential forms enables one to describe discrete 
transitions, evolutionary processes and generation of various structures.

\section{Identical relations and nondegenerate transformations}

Identical relations lie at the basis of the mathematical 
apparatus of exterior differential forms. They reflect the properties 
of exterior forms. Below we present some properties of closed exterior 
differential forms, which are necessary for further presentation. 
(In more detail about skew-symmetric differential forms one can read 
in [1-7]).

\subsection*{Closed exterior differential forms}

The exterior differential form of degree $p$ ($p$-form on the 
differentiable manifold) can be written as [1,3] 
$$
\theta^p=\sum_{i_1\dots i_p}a_{i_1\dots i_p}dx^{i_1}\wedge
dx^{i_2}\wedge\dots \wedge dx^{i_p}\quad 0\leq p\leq n\eqno(1.1)
$$
Here $a_{i_1\dots i_p}$ are the functions of the variables $x^{i_1}$, 
$x^{i_2}$, \dots, $x^{i_p}$, $n$ is the dimension of space, 
$\wedge$ is the operator of exterior multiplication, $dx^i$, 
$dx^{i}\wedge dx^{j}$, $dx^{i}\wedge dx^{j}\wedge dx^{k}$, \dots\
is the local basis, which satisfies the condition of exterior 
multiplication:
$$
\begin{array}{l}
dx^{i}\wedge dx^{i}=0\\
dx^{i}\wedge dx^{j}=-dx^{j}\wedge dx^{i}\quad i\ne j
\end{array}\eqno(1.2)
$$

The differential of the (exterior) form $\theta^p$ is expressed as 
$$
d\theta^p=\sum_{i_1\dots i_p}da_{i_1\dots
i_p}dx^{i_1}dx^{i_2}\dots dx^{i_p} \eqno(1.3)
$$
and is the differential form of degree $(p+1)$. 

[From here on the symbol $\sum$ can be omitted and it will be 
implied that a summation over double indices is performed.  Besides, the 
symbol of exterior multiplication will be also omitted for the 
sake of presentation convenience]. 

Let us consider the exterior differential form of the first degree 
$\omega=a_i dx^i$. In this case the differential will be expressed
as $d\omega=K_{ij}dx^i dx^j$, where 
$K_{ij}=(\partial a_j/\partial x^i-\partial a_i/\partial x^j)$ are 
components of the form commutator. 

In this section we will consider the domains of Euclidean space or 
differentiable manifolds [2]. (Manifolds, on 
which the skew-symmetric differential forms may be 
defined, and the influence of the manifold properties on the 
differential forms  will be discussed in more detail in section 2). 

The exterior differential form of degree $p$ ($p$-form on the 
differentiable manifold) is called a closed one 
if its differential is equal to zero: 
$$
d\theta^p=0\eqno(1.4)
$$

A differential of the form is a closed form. That is, 
$$
dd\omega=0\eqno(1.5)
$$
where $\omega$ is an arbitrary exterior form. 

The form that is a differential of some other form: 
$$
\theta^p=d\theta^{p-1}\eqno(1.6)
$$
is called an {\it exact} form. Exact forms prove to be closed 
automatically [2] 
$$
d\theta^p=dd\theta^{p-1}=0\eqno(1.7)
$$

Here it is necessary to call attention to the following points. In the 
above presented formulas it was implicitly assumed that the 
differential operator $d$ is a total one (that is, the 
operator $d$ acts everywhere in the vicinity of the point 
considered locally),  and therefore, it acts on the manifold of the 
initial dimension $n$. However, a differential may be internal. Such a 
differential acts on some structure with the dimension being less 
than that of the initial manifold. The structure, on which the exterior 
differential form may become a closed {\it inexact} form, is a 
pseudostructure with respect to its metric properties. \{Cohomology, 
sections of cotangent bundles, the eikonal surfaces, the 
characteristical and potential surfaces, 
and so on may be regarded as examples of pseudostructures.\} 

If the form is closed on  pseudostructure only, the closure 
condition is written as 
$$
d_\pi\theta^p=0\eqno(1.8)
$$
And the pseudostructure $\pi$ is defined from the condition 
$$
d_\pi{}^*\theta^p=0\eqno(1.9)
$$
where ${}^*\theta^p$ is a dual form. 
(For the properties of dual forms see [5]). 

The fundamental properties of exterior differential forms are connected 
with the fact that any closed form is a differential. 
The exact form is, by definition, a differential (see condition (1.6)). 
In this case the differential is total. The closed inexact form is 
a differential too. The closed inexact form is an interior (on 
pseudostructure) differential, that is 
$$
\theta^p_\pi=d_\pi\theta^{p-1}\eqno(1.10)
$$

And so, any closed form is a differential of the form of lower 
degree: the total one $\theta^p=d\theta^{p-1}$ if the form is exact, 
or the interior one $\theta^p=d_\pi\theta^{p-1}$ on pseudostructure if 
the form is inexact. 

The closure of exterior differential forms result from the conjugacy of 
elements of exterior or dual forms. 

The closure property of the exterior form means that any objects, namely, 
elements of the exterior form, components of elements, elements of 
the form differential, exterior and dual forms and others, turn 
out to be conjugated. A  variety of objects of conjugacy leads to the 
fact that there is a large number of different types of closed 
exterior forms. 

Since the conjugacy is a certain connection between two operators or 
mathematical objects, it is evident that, to express a conjugacy 
mathematically, it can be used relations. Just such relations constitute 
the basis of mathematical apparatus of the exterior differential forms. 
This is an identical relation.

\subsection*{Identical relations of exterior differential forms} 

The identical relations reflect the closure conditions 
of the differential forms, namely, vanishing the form 
differential (see formulas (1.4), (1.8), (1.9)) and the conditions 
connecting the forms of consequent degrees (see formulas (1.6), (1.10)). 

The importance of the identical relations for exterior differential 
forms is manifested by the fact that practically in all branches of 
physics, mechanics, thermodynamics one faces such identical relations. 
One can present the following examples: 

a) the Poincare invariant $ds\,=\,-H\,dt\,+\,p_j\,dq_j$, 

b) the second principle of thermodynamics $dS\,=\,(dE+p\,dV)/T$, 

c) the vital force theorem in theoretical mechanics: $dT=X_idx^i$ 
where $X_i$ are the components of potential force, and $T=mV^2/2$ is 
the vital force, 

d) the conditions on characteristics in the theory of differential 
equations, and so on. 

The identical relations in differential forms express the fact that each 
closed exterior form is a differential of some exterior form (with the 
degree less by one). In general form such an identical relation can be 
written as 
$$
d _{\pi}\phi=\theta _{\pi}^p\eqno(1.11)
$$
In this relation the form in the right-hand side has to be a {\it closed} 
one. (As it will be shown below, the identical relations are satisfied 
only on pseudostructures). 

In identical relation (1.11) in one side it stands the closed form and 
in other side does the  differential of some differential form of 
the less by one degree, which is the closed form as well. 

In addition to relations in the differential forms from the closure 
conditions of differential forms and the conditions connecting 
the forms of consequent degrees  the identical relations of other 
types are obtained. The types of such relations are presented below. 

1. {\it Integral identical relations}. 

The formulas by Newton, Leibnitz, Green, the integral relations by 
Stokes, Gauss-Ostrogradskii are examples of integral identical relations. 

2. {\it Tensor identical relations}. 

From the relations that connect exterior forms of consequent degrees 
one can obtain the vector and tensor identical relations that connect 
the operators of the gradient, curl, divergence and so on. 

From the closure conditions of exterior and dual forms one can obtain 
the identical relations such as the gauge relations in electromagnetic 
field theory, the tensor relations between connectednesses and their 
derivatives in gravitation (the symmetry of connectednesses with respect 
to lower indices, the Bianchi identity, the conditions imposed on the 
Christoffel symbols) and so on.

3. {\it Identical relations between derivatives}. 

The identical relations between derivatives correspond to the closure 
conditions of exterior and dual forms. The examples of such relations 
are the above presented Cauchy-Riemann conditions in the theory of 
complex variables, the transversality condition in the calculus of 
variations, the canonical relations in the Hamilton 
formalism, the thermodynamic relations between derivatives 
of thermodynamic functions [8], the condition which the derivative of 
implicit function is subject to, the eikonal relations [9] and so on.

\subsection*{Nondegenerate transformations} 

One of the fundamental methods in the theory of exterior 
differential forms is an application of {\it nondegenerate} 
transformations (below it will be said about {\it degenerate} 
transformations). 

In the theory of exterior differential 
forms the nondegenerate transformations are those that conserve the 
differential. This is connected with the property of closed 
differential forms. Since a closed form is a differential (a total one, 
if the form is exact, or an interior one on pseudostructure, if 
the form is inexact), it is evident that the closed form turns out to 
be invariant under all transformations that conserve the differential.  

The examples of nondegenerate transformations in the theory  of 
exterior differential forms are unitary, tangent, canonical, gradient 
transformations.

To the nondegenerate transformations there are assigned closed forms 
of given degree. To the unitary transformations it is assigned (0-form), 
to the tangent and canonical transformations it is assigned (1-form), to 
the gradient transformations it is assigned (2-form) and so on. 
It should be noted that these transformations are {\it gauge 
transformations} for spinor, scalar, vector, tensor (3-form) fields. 

The connection between nondegenerate transformations and closed exterior 
forms disclose an internal commonness of  nondegenerate transformations: 
all these  transformations are transformations that preserve 
a differential and nondegenerate transformations can be classified by 
a degree of corresponding closed differential or dial forms. 

Nondegenerate transformations, if applied to identical relations, 
enable one to obtain new identical relations and new closed exterior 
differential forms.

The nondegenerate transformations proceed the transitions between 
the closed forms of different degrees. With such nondegenerate 
transformations of the exterior differential forms many operators 
of mathematical physics are connected. If, in addition to the exterior 
differential, we introduce the following operators: 1) $\delta$ for 
transformations that convert the form of $(p+1)$ degree into the form 
of $p$ degree, 2) $\delta'$ for cotangent transformations, 3) $\Delta$ 
for the $d\delta-\delta d$ transformation, 4) $\Delta'$ for the 
$d\delta'-\delta'd$ transformations, one can see that the operator 
$\delta$ corresponds to Green's operator, $\delta'$ does to the 
canonical transformation operator, $\Delta$ does to the d'Alembert 
operator in 4-dimensional space, and $\Delta'$ corresponds to the 
Laplace operator [5]. 

Hence, one can see that many identical relations and nondegenerate 
transformations, which occur in various branches of mathematics, are 
connected with the properties of closed exterior (skew-symmetric) 
differential forms. This enables one to understand the properties and  
specific features of such identical relations and nondegenerate 
transformations and to introduce their classification. (And yet 
this discloses potentialities of exterior differential forms.) 

The role of identical relations in mathematics and physics 
consists in the fact that they define invariant and covariant 
objects, which are of great functional and utilitarian importance. 
The above considered identical relations, which correspond to 
closed exterior and dual forms, define the objects, which remain 
to be invariant and covariant under all transformations preserving 
the differential. I should be noted that an invariant object, to 
which the closed form corresponds, and the covariant object, to 
which the dual form corresponds, form a double structure being the 
example of the differential and geometric G-structure. In physics 
such structures describe the physical structures that forms 
physical fields and relevant manifolds [10,11].

\section{Nonidentical relations and degenerate transformations}

Identical relations and nondegenerate transformations, which describe 
invariant objects and transitions between invariant objects, cannot 
explain how do the invariant objects appear and what do them generate. 
To answer this question one must understand how do the identical 
relations realize. The theory of exterior differential forms, 
which is invariant one, cannot give the answer to such questions. 
To do this, the evolutionary theory is necessary. In the author's 
works it has been sown that there exist a skew-symmetric differential 
forms, which possess the evolutionary properties. They form nonidentical 
relations from which the identical relations, corresponding to invariant 
objects, are obtained by application of the degenerate transformations.

\subsection*{Evolutionary differential forms} 
Differential forms, which possess the evolutionary properties and so 
were called the evolutionary differential forms, appear 
from the description of various process. 

A radical  distinction between the evolutionary forms and the exterior 
ones consists in the fact that the exterior differential forms are 
defined on manifolds with {\it closed metric forms}, whereas the 
evolutionary differential forms are defined on manifolds with 
{\it unclosed metric forms}. 

The evolutionary differential form of degree $p$ ($p$-form), 
as well as the exterior differential form, can be written down as 
$$
\omega^p=\sum_{\alpha_1\dots\alpha_p}a_{\alpha_1\dots\alpha_p}dx^{\alpha_1}\wedge
dx^{\alpha_2}\wedge\dots \wedge dx^{\alpha_p}\quad 0\leq p\leq n\eqno(2.1)
$$
But the evolutionary form differential cannot be written similarly to 
that presented for exterior differential forms (see formula (1.3)). 
In the evolutionary form differential there appears an additional term 
connected with the fact that the basis of the form changes. For 
differential forms defined on the manifold with unclosed metric form one 
has $d(dx^{\alpha_1}dx^{\alpha_2}\dots dx^{\alpha_p})\neq 0$ (it should 
be noted that for differentiable manifold the following is valid:
$d(dx^{\alpha_1}dx^{\alpha_2}\dots dx^{\alpha_p}) = 0$). 
For this reason the differential of the evolutionary form $\omega^p$ 
can be written as
$$
d\omega^p{=}\!\sum_{\alpha_1\dots\alpha_p}\!da_{\alpha_1\dots\alpha_p}dx^{\alpha_1}dx^{\alpha_2}\dots
dx^{\alpha_p}{+}\!\sum_{\alpha_1\dots\alpha_p}\!a_{\alpha_1\dots\alpha_p}d(dx^{\alpha_1}dx^{\alpha_2}\dots
dx^{\alpha_p})\eqno(2.2)
$$
where the second term is connected with the differential of the basis. 
The second term is expressed in terms of the metric form commutator. 
For the manifold with closed metric form this term vanishes.

For example, we again inspect the first-degree form 
$\omega=a_\alpha dx^\alpha$. The differential of this form can 
be written as $d\omega=K_{\alpha\beta}dx^\alpha dx^\beta$, where 
$K_{\alpha\beta}=a_{\beta;\alpha}-a_{\alpha;\beta}$ are 
components of the commutator of the form $\omega$, and 
$a_{\beta;\alpha}$, $a_{\alpha;\beta}$ are the covariant 
derivatives. If we express the covariant derivatives in terms of 
the connectedness (if it is possible), they can be written 
as $a_{\beta;\alpha}=\partial a_\beta/\partial 
x^\alpha+\Gamma^\sigma_{\beta\alpha}a_\sigma$, where the first 
term results from differentiating the form coefficients, and the 
second term results from differentiating the basis. (In 
Euclidean space covariant derivatives coincide with ordinary ones 
since in this case derivatives of the basis vanish). If we substitute 
the expressions for covariant derivatives into the formula for 
the commutator components, we obtain the following expression for 
the commutator components of the evolutionary form $\omega$: 
$$
K_{\alpha\beta}=\left(\frac{\partial a_\beta}{\partial
x^\alpha}-\frac{\partial a_\alpha}{\partial
x^\beta}\right)+(\Gamma^\sigma_{\beta\alpha}-
\Gamma^\sigma_{\alpha\beta})a_\sigma\eqno(2.3)
$$
Here the expressions 
$(\Gamma^\sigma_{\beta\alpha}-\Gamma^\sigma_{\alpha\beta})$ 
entered into the second term are just the components of 
commutator of the first-degree metric form. 

That is, the corresponding metric form commutator 
will enter into the differential form commutator. 

The evolutionary properties of evolutionary differential forms are 
just connected with the properties of commutator of this form.  

The evolutionary differential form commutator, in contrast to that of 
the exterior one, cannot be equal to zero because it involves the metric 
form commutator being nonzero. This means that the evolutionary form 
differential is nonzero. Hence, the evolutionary differential form, in 
contrast to the case of the exterior form,  cannot be closed. 

'he commutators of evolutionary forms depend not only on the 
evolutionary form coefficients, but on the characteristics of manifolds, 
on which this form is defined, as well. As a result, such a dependence 
of the evolutionary form commutator produces the topological and 
evolutionary properties of both the commutator and the evolutionary 
form itself (this will be demonstrated below). 

Since the evolutionary differential forms are unclosed, the mathematical 
apparatus of evolutionary differential forms does not seem to 
possess any possibilities connected with the algebraic, group, invariant 
and other properties of closed exterior differential forms. However, the 
mathematical apparatus of evolutionary forms includes some new 
unconventional elements.

\subsection*{Nonidentical relations of evolutionary differential forms} 

Above it was shown that the identical relations lie at the 
basis of the mathematical apparatus of exterior differential forms. 

In contrast to this, nonidentical relations lie at the basis of the 
mathematical apparatus of evolutionary differential forms. 

The identical relations of closed exterior differential forms reflect 
a conjugacy of any objects. The evolutionary forms, being unclosed, 
cannot directly describe a conjugacy of any objects. But they allow 
a description of the process in which the conjugacy may appear (the 
process when closed exterior differential forms are generated). Such 
a process is described by nonidentical relations. 

The concept of ``nonidentical relation"  may appear to be inconsistent. 
However, it has a deep meaning. 

The identical relations establish exact correspondence between the 
quantities (or objects) involved into this relation. This is possible 
in the case when the quantities involved into the relation are 
measurable ones. [A quantity is called a measurable quantity if its 
value does not change under transition to another, equivalent, 
coordinate system. In other words, this quantity is invariant one.] 
In the nonidentical relations one of the quantities is unmeasurable. 
(Nonidentical relations with two unmeasurable quantities are senseless). 
If this relation is evolutionary one, it turns out to be a selfvarying 
relation, namely, a variation of some object leads to  variation of 
other one, and in turn a variation of the second object leads to 
variation of the first and so on. Since in the nonidentical relation 
one of the objects is a unmeasurable quantity, the other object cannot 
be compared with the first, and therefore, the process cannot stop. 
Here the specific feature is that in the process of such selfvarying 
it may be realized the additional conditions under which the identical 
relation can be obtained from nonidentical relation. 
The additional condition may be realized spontaneously while 
selfvarying the nonidentical relation if the system (which is described 
by such relation) possesses any symmetry. When such additional 
conditions are realized the exact correspondence between the quantities 
involved in the relation is established. Under the additional condition 
a unmeasurable quantity becomes a measurable quantity as well, and
the exact correspondence between the objects involved in the 
relation is established. That is, an identical relation 
can be obtained from a nonidentical relation. 

The nonidentical relation is a relation between a closed exterior 
differential form, which is a differential and is a measurable quantity, 
and an evolutionary form, which is an unmeasurable quantity. 

Nonidentical relations of such type appear under descriptions of any 
processes. These relations may be written as 
$$
d\psi \,=\,\omega^p \eqno(2.4)
$$
Here $\omega^p$ is the $p$-degree evolutionary form that is 
nonintegrable, $\psi$ is some form of degree $(p-1)$, and 
the differential $d\psi$ is a closed form of degree $p$. 

In the left-hand side of this relation it stands the form differential, 
i.e. the closed form, which is an invariant object. In the right-hand 
side it stands the nonintegrable unclosed form, which is not an 
invariant object. Such a relation cannot be identical. 

One can see a difference of relations for exterior forms and 
evolutionary ones. 
In the right-hand side of identical relation (see relation (1.11)) 
it stands a closed form, whereas the form in the right-hand side of 
nonidentical relation (2.4) is an unclosed one. 

Such nonidentical relations appear, for example, while 
investigating the integrability of any differential equations that 
describe various processes. The equation is integrable if it can 
be reduced to the form $d\psi=dU$. However it appears that, if the 
equation is not subject to an additional condition (the 
integrability condition), it is reduced to form (2.4), where 
$\omega$ is an unclosed form and it cannot be expressed as a 
differential. In the works [6,10] there are presented nonidentical 
relations obtained while investigating the integrability of the 
first-order partial differential equation and while analyzing the 
physical processes in material media. (it could be noted that the 
differential in the left-hand side of the nonidentical relation 
specifies the state of the medium or system described, and $d\psi$ 
can be a state function.) 

The first principle of thermodynamics is an example of nonidentical 
relation [8]. 

{\bf\it  It arises a question of how to work with nonidentical relation?} 

Two different approaches are possible. 

The first, evident, approach is to find the condition under which the 
nonidentical relation becomes identical and to obtain a closed form 
under this condition. In other words, the nonidentical relation 
is subject to the condition, under which this relation is transformed 
into an identical relation (if it is possible). 

Such an approach is traditional and is always used implicitly. It may 
be shown that additional conditions are imposed on the mathematical 
physics equations obtained in description of the physical processes so 
that these equations should be invariant (integrable) or should have 
invariant solutions. 

This approach does not solve the evolutionary problem. 
\{Here a psychological point should be noted. While investigating real 
physical processes one often faces the relations that are nonidentical. 
But it is commonly believed that only identical relations can have any 
physical meaning. For this reason one immediately attempts to impose 
a condition onto the nonidentical relation under which this relation 
becomes identical, and it is considered only such cases when this 
relation can satisfy the additional conditions. And all remaining is 
rejected. It is not taken into account that a nonidentical relation is 
often obtained from a description of some physical process and it has 
physical meaning at every stage of the physical process rather 
than at the stage when the additional conditions are satisfied.  
In essence, the physical process does not considered completely. 
At this point it should be emphasized that the nonidentity of 
the evolutionary relation does not mean the imperfect accuracy 
of the mathematical description of a physical process. The 
nonidentical relations are indicative of specific features of 
the physical process development.\} 

The second approach to investigating the nonidentical evolutionary 
relation shows that the conditions, under which from the nonidentical 
relation it is obtained the identical relation, are realized under 
selfvariation of the nonidentical relation.

\subsection*{Selfvariation of the evolutionary nonidentical relation} 
The evolutionary nonidentical relation is selfvarying, 
because, firstly, it is nonidentical, namely, it contains 
two objects one of which appears to be unmeasurable, and, 
secondly, it is an evolutionary relation, namely, a variation of 
any object of the relation in some process leads to variation of 
another object and, in turn, a variation of the latter leads to 
variation of the former. Since one of the objects is an unmeasurable 
quantity, the other cannot be compared with the first one, and hence, 
the process of mutual variation cannot stop.

Varying the evolutionary form coefficients leads to varying the first 
term of the evolutionary form commutator (see (2.3)). In accordance with 
this variation it varies the second term, that is, the metric form of 
manifold varies. Since the metric form commutators specifies the 
manifold differential characteristics, which are connected with the 
manifold deformation (as it has been pointed out, the commutator of the 
zero degree metric form specifies the bend, that of second degree 
specifies various types of rotation, that of the third degree specifies 
the curvature), this points to the manifold deformation. This means that 
it varies the evolutionary form basis. In turn, this leads to variation 
of the evolutionary form, and the process of intervariation of the 
evolutionary form and the basis is repeated. Processes of variation of 
the evolutionary form and the basis are controlled by the evolutionary 
form commutator and it is realized according to the evolutionary relation. 

The process of the evolutionary relation selfvariation plays a governing 
role in description of the evolutionary processes. 

The significance of the evolutionary relation selfvariation consists in 
the fact that in such a process it can be realized the conditions under 
which the identical relation is obtained from the nonidentical relation. 
These are the conditions of degenerate transformation.

\subsection*{Degenerate transformations.} 

To obtain the identical relation from the evolutionary nonidentical 
relation, it is necessary that a closed exterior differential form 
should be derived from the evolutionary differential form, which is 
included into evolutionary  relation. However, as it has been shown 
above, the evolutionary form cannot be a closed form. For this reason 
a transition from the evolutionary form is possible only to an 
{\it inexact} closed exterior form, which is defined on pseudostructure. 

To the pseudostructure it is assigned a closed dual form 
(whose differential vanishes). For this reason a transition 
from the evolutionary form to a closed inexact exterior form proceeds 
only when the conditions of vanishing the dual form differential are 
realized, in other words, when the metric form differential or 
commutator becomes equal to zero. 

Since the evolutionary form differential is nonzero, whereas the closed 
exterior form differential is zero, the transition from the evolutionary 
form to the closed exterior form is allowed only under {\it degenerate 
transformation}. The conditions of vanishing the dual form differential 
(the additional condition) are the conditions of degenerate 
transformation. 

Such conditions can just be realized under selfvariation of the 
nonidentical evolutionary relation. 

As the conditions of degenerate transformation (additional conditions) 
it can serve any symmetries of the evolutionary form coefficients 
or of its commutator. (While describing material system such additional 
conditions are related, for example, to 
degrees of freedom of the material system). 

Mathematically to the conditions of degenerate transformation there 
corresponds a requirement that some functional expressions become equal 
to zero. Such functional expressions are Jacobians, determinants, 
the Poisson brackets, residues, and others.

\subsection*{Obtaining identical relation from  nonidentical one} 

Let us consider nonidentical evolutionary relation (2.4). 

As it has been already mentioned, the evolutionary differential form 
$\omega^p$, involved into this relation is an unclosed one. The 
commutator, and hence the differential, of this form is nonzero. That is, 
$$
d\omega^p\ne 0\eqno(2.5)
$$
If the transformation is degenerate, from the unclosed evolutionary form 
it can be obtained a differential form closed on pseudostructure. 
The differential of this form equals zero. That is, it is 
realized the transition 

 $d\omega^p\ne 0 \to $ (degenerate transformation) $\to d_\pi \omega^p=0$, 
$d_\pi{}^*\omega^p=0$ 

On the pseudostructure $\pi$ evolutionary relation (2.4) transforms into 
the relation 
$$
d_\pi\psi=\omega_\pi^p\eqno(2.6)
$$
which proves to be the identical relation. Indeed, since the form 
$\omega_\pi^p$ is a closed one, on the pseudostructure this form turns 
out to be a differential of some differential form. In other words, 
this form can be written as $\omega_\pi^p=d_\pi\theta$. Relation (2.11) 
is now written as 
$$
d_\pi\psi=d_\pi\theta
$$
There are differentials in the left-hand and right-hand sides of 
this relation. This means that the relation is an identical one. 

From evolutionary relation (2.4) it is obtained the identical on the 
pseudostructure relation. In this case the evolutionary relation itself 
remains to be nonidentical one. (At this point it should be 
emphasized that differential, which equals zero, is an interior one. 
The evolutionary form commutator becomes zero only on the 
pseudostructure. The total evolutionary form commutator is nonzero. That 
is, under degenerate transformation the evolutionary form differential 
vanishes only {\it on pseudostructure}. The total differential of the 
evolutionary form is nonzero. The evolutionary form remains to be 
unclosed.) 

It can be shown that all identical relations of the exterior 
differential form theory are obtained from nonidentical relations (that 
contain the evolutionary forms) by applying degenerate transformations. 

{\it The degenerate transform is realized as a transition to 
nonequivalent coordinate system: the transition from the accompanying 
noninertial coordinate system to the locally inertial that}. 
Evolutionary relation (2.4) and condition (2.5) relate 
to the system being tied to the accompanying manifold, whereas 
identical relations (2.6) may relate only to the locally inertial 
coordinate system being tied to a pseudostructure. 

\subsection*{Integration of the nonidentical evolutionary relation} 
Under degenerate transformation from the nonidentical evolutionary 
relation one obtains a relation being identical on pseudostructure. 
Since the right-hand side of such a relation can be expressed in terms 
of differential (as well as the left-hand side), one obtains a relation 
that can be integrated, and as a result he obtains a relation with the 
differential forms of less by one degree. 

The relation obtained after integration proves to be nonidentical as 
well. 

The resulting nonidentical relation of degree $(p-1)$ (relation that 
contains the forms of the degree $(p-1)$) can be integrated once again 
if the corresponding degenerate transformation has been realized and 
the identical relation has been formed. 

By sequential integrating the evolutionary relation of degree $p$ (in 
the case of realization of the corresponding degenerate transformations 
and forming the identical relation), one can get closed (on the 
pseudostructure) exterior forms of degree $k$, where $k$ ranges from 
$p$ to $0$. 

In this case one can see that under such integration the closed (on the 
pseudostructure) exterior forms, which depend on two parameters, are 
obtained. These parameters are the degree of evolutionary form $p$ (in 
the evolutionary relation) and the degree of created closed forms $k$. 

In addition to these parameters, another parameter appears, namely, the 
dimension of space. If the evolutionary relation generates the closed 
forms of degrees $k=p$, $k=p-1$, \dots, $k=0$, to them there correspond 
the pseudostructures of dimensions $(N-k)$, where $N$ is the space 
dimension. \{It is known that to the closed exterior differential forms 
of degree $k$ there correspond skew-symmetric tensors of rank $k$ and to 
corresponding dual forms there do the pseudotensors of rank $(N-k)$, 
where $N$ is the space dimensionality. The pseudostructures correspond 
to such tensors, but only on the space formed.\}

\subsection*{Relation between degenerate and nondegenerate transformations} 

A peculiarity of the degenerate and nondegenerate transformations can be 
considered by the example of the field equation (the investigation of 
this equation is presented in the work [6]). In this case the degenerate 
transformation is a transition from the 
Lagrange function to the Hamilton function. The equation for the 
Lagrange function, that is the Euler variational equation, was obtained 
from the condition $\delta S\,=\,0$, where $S$ is the action functional. 
In the real case, when forces are nonpotential or couplings are 
nonholonomic, the quantity $\delta S$ is not a closed form, that is, 
$d\,\delta S\,\neq \,0$. But the Hamilton function is obtained from 
the condition $d\,\delta S\,=\,0$ which is the closure condition for 
the form $\delta S$. The transition from the Lagrange function $L$ to 
the Hamilton function $H$ (the transition from variables $q_j,\,\dot q_j$ 
to variables $q_j,\,p_j=\partial L/\partial \dot q_j$) is a transition 
from the tangent space, where the form is unclosed, to the cotangent 
space with a closed form. One can see that this transition is 
a degenerate one. 

The degenerate transformation is a transition from the tangent space 
($q_j,\,\dot q_j)$) to the cotangent (characteristic) manifold 
($q_j,\,p_j$). On the other hand, the nondegenerate canonical 
transformation is a transition from one characteristic manifold 
($q_j,\,p_j$) to another characteristic manifold ($Q_j,\,P_j$). 
$\{$The formula of nondegenerate canonical 
transformation can be written as $p_jdq_j=P_jdQ_j+dW$, 
where $W$ is the generating function$\}$. 

Here it has been shown a connection between the canonical nondegenerate 
transformation and the degenerate transformation. 
It may be easily shown that such a property of duality is also a 
specific feature of transformations such as tangent, gradient, 
contact, gauge, conform mapping, and others.

\subsection*{Transition from nonconjugated operators to conjugated 
operators} 

The evolutionary process of obtaining the identical relation from 
the nonidentical one and obtaining a closed (inexact) exterior 
form from the unclosed evolutionary form describes a process of 
conjugating any objects. 

In section 1 it has been shown that the condition of the closure 
of exterior differential forms is a result of the conjugacy of 
any constituents  of the exterior or dual forms (the form elements, 
components of each element, exterior and dual forms, exterior forms of 
various degrees, and others). Since the identical 
relations of exterior differential forms is a mathematical  
record of the closure conditions of exterior differential forms and, 
correspondingly, of conjugacy of any objects, the process of obtaining 
the identical relation from nonidentical one (selfmodification of the 
nonidentical evolutionary relation and degenerate transformation) is a 
process of conjugating any objects. 

Here it could be pointed out the following. 
To the differential of the  closed exterior differential form 
there correspond conjugated operators being equal to zero, whereas to 
the differential of the evolutionary form there correspond 
nonconjugated operators being not equal to zero. The transition from 
the evolutionary form to the closed exterior form is that from 
nonconjugated operators to conjugated ones. This is expressed 
mathematically as a transition from a nonzero differential 
(the evolutionary form differential is nonzero) to a differential that 
equals zero (the closed exterior form differential equals zero). 
'his reveals as the transition from one coordinate system to another 
(nonequivalent) coordinate system.

\bigskip
{\large\bf Summary}.

The mathematical apparatus of exterior and evolutionary skew-symmetric 
differential forms constitute a new closed mathematical apparatus that 
possesses the unique properties. It includes new, unconventional, 
elements: "nonidentical relation", "degenerate transformation", 
"transition from one frame of reference to another, nonequivalent, 
frame of reference". This enables one to create the mathematical 
language that has radically new abilities.

Identical and nonidentical relations, nondegenerate  and 
degenerate transformations, transitions from nonidentical 
relations to identical ones are of great functional and 
utilitarian importance. First of all, this importance consists in 
the fact that, using the degenerate transformations, the 
nonidentical relations generate the identical relations and closed 
exterior forms, which describe invariant and covariant objects and 
the differential and geometrical structures. This discloses the 
mechanism of evolutionary processes, discrete transitions, 
generation of various structures. 

The utilitarian importance of nonidentical relations and degenerate 
transformations consists in the fact that they disclose a mechanism 
of the evolutionary processes in material media and a process of 
creating physical structures and forming physical fields and manifolds. 
In more detail this is outlined in the works [10,11].

1. Bott R., Tu L.~W., Differential Forms in Algebraic Topology. 
Springer, NY, 1982. 

2. Schutz B.~F., Geometrical Methods of Mathematical Physics. Cambrige 
University Press, Cambrige, 1982. 

3. Encyclopedia of Mathematics. -Moscow, Sov.~Encyc., 1979 (in Russian). 

4. Novikov S.~P., Fomenko A.~P., Elements of the differential geometry 
and topology. -Moscow, Nauka, 1987 (in Russian). 

5. Wheeler J.~A., Neutrino, Gravitation and Geometry. Bologna, 1960. 

6. Petrova L.~I., Exterior and evolutionary skew-symmetric differential 
forms and their role in mathematical physics. 

http://arXiv.org/pdf/math-ph/0310050 

7. Petrova L.~I., Invariant and evolutionary properties of the 
skew-symmetric differential forms. 

http://arXiv.org/pdf/math.GM/0401039

8. Haywood R.~W., Equilibrium Thermodynamics. Wiley Inc. 1980. 

9. Fock V.~A., Theory of space, time, and gravitation. -Moscow, 
Tech.~Theor.~Lit., 1955 (in Russian).

10. Petrova L.~I., Conservation laws. Their role in evolutionary 
processes. (The method of skew-symmetric differential forms).

http://arXiv.org/pdf/math-ph/0311008 

11. Petrova L.~I., Approaches to general field theory. 
(The method of skew-symmetric differential forms) 

http://arXiv.org/pdf/math-ph/0311030 

\end{document}